\documentclass[12pt]{article}
\usepackage{epstopdf}
\usepackage{graphicx}
\usepackage[ruled,vlined]{algorithm2e}
\usepackage{appendix}
\usepackage{amsmath}
\usepackage{amssymb}
\usepackage{amsthm}
\usepackage{multirow}
\usepackage{color}
\usepackage{longtable}
\usepackage{array}
\usepackage{url}
\usepackage{booktabs}
\usepackage{float}
\usepackage{mathtools}
\usepackage{tikz}
\usepackage{caption}
\usepackage{graphicx}
\usepackage{float}
\usepackage{subfigure}
\usepackage{bm}

\newtheorem{theorem}{Theorem}[section]

\theoremstyle{definition}

\title{A note on multicolor Ramsey number of small odd cycles versus a large clique}
\author{Zixiang Xu$^{\text{a,}}$\thanks{e-mail: zxxu8023@qq.com.} and Gennian Ge$^{\text{a,}}$\thanks{e-mail: gnge@zju.edu.cn. Research supported by the National Key Research and Development Program of China under Grant No. 2020YFA0712100, National Natural Science Foundation of China under Grant No. 11971325, and Beijing Scholars Program.}\\
\footnotesize $^{\text{a}}$ School of Mathematical Sciences, Capital Normal University, Beijing 100048, China.
}

\begin{document}

\date{}

\maketitle

\begin{abstract}
 Let $R_k(H;K_m)$ be the smallest number $N$ such that every coloring of the edges of $K_{N}$ with $k+1$ colors has either a monochromatic $H$ in color $i$ for some $1\leqslant i\leqslant k$, or a monochromatic $K_{m}$ in color $k+1$. In this short note, we study the lower bound for $R_k(H;K_m)$ when $H$ is $C_5$ or $C_7$, respectively. We show that
  \begin{equation*}
      R_{k}(C_5;K_m)=\Omega(m^{\frac{3k}{8}+1}/(\log{m})^{\frac{3k}{8}+1}),
  \end{equation*} and
  \begin{equation*}
      R_{k}(C_7;K_m)=\Omega(m^{\frac{2k}{9}+1}/(\log{m})^{\frac{2k}{9}+1}),
  \end{equation*}
  for fixed positive integer $k$ and $m\rightarrow\infty$.
  These slightly improve the previously known lower bound
  $R_{k}(C_{2\ell+1};K_m)=\Omega(m^{\frac{k}{2\ell-1}+1}/(\log m)^{k+\frac{2k}{2\ell-1}})$ obtained by Alon and R\"{o}dl. The proof is based on random block constructions and random blowups argument.


\smallskip

\end{abstract}

\section{Introduction}
For a given graph $H$ and a positive integer $k$, the multicolor Ramsey number $R_k(H;K_m)$ is the smallest number $N$ such that every edge coloring of $K_{N}$ with $k+1$ colors has either a monochromatic $H$ in color $i$ for some $1\leqslant i\leqslant k$, or a monochromatic $K_{m}$ in color $k+1$. Obtaining good upper and lower bounds on $R_k(H;K_m)$ is a central question in Ramsey theory,
and despite more than $80$ years of effort, still not much is known. For example, the most well-studied case
is that of $k=1$ and $H=K_{m}$, and the theoretical bounds $2^{\frac{m}{2}}\leqslant R_{1}(K_{m};K_{m})\leqslant 2^{2m}$ were obtained by Erd\H{o}s and Szekeres~\cite{1947Erdos, 1935Erdos} in $1940s$. However, only lower term improvements~\cite{2009Conlon, 2020Sah, 1975Joel} have been made to these bounds. For larger values of $k\geqslant 2,$ even less has been known. When $H=K_{m}$, a recent breakthrough of Conlon and Ferber~\cite{2021conlon} yielded an exponential improvement on the lower bounds for multicolor diagonal Ramsey numbers $R_{k}(K_{m};K_{m})$ using mixed random algebraic approaches. More recently, the lower bounds for $R_{k}(K_{m};K_{m})$ were improved by Wigderson~\cite{2021Yuval} and Sawin~\cite{2021Sawin} via different methods.

In this paper, we mainly study the lower bound for $R_{k}(H;K_{m})$ when $H$ is an odd cycle $C_{2\ell+1}$. For $k=1$, Erd\H{o}s, Faudree, Rousseau and Schelp~\cite{1978Erdos} proved that $R_{1}(C_{\ell};K_{m})=O(m^{1+\frac{1}{t}})$, where $t=\lceil\frac{\ell}{2}\rceil-1$. These were later improved by Sudakov~\cite{2002Sudakov}, who proved that $R_{1}(C_{2\ell+1};K_{m})=O(\frac{m^{1+1/\ell}}{\log^{1/\ell}m})$. When $\ell=1$, the matched lower bound for $R_{1}(C_{3},K_{m})=\Theta(\frac{m^{2}}{\log{m}})$ was shown in~\cite{1980Ajtai, 1995Kim}. Moreover, by analyzing the $C_{\ell}$-free process, Bohman and Keevash~\cite{2010Keevash} proved that $R_{1}(C_{\ell};K_{m})=\Omega(\frac{m^{(\ell-1)/(\ell-2)}}{\log{m}})$. For $\ell\in\{5,6,7,10\}$, the best known lower bounds for $R_{1}(C_{\ell};K_{m})$ were improved by Mubayi and Verstra\"{e}te~\cite{2019Mubayi} via different approaches. When $k\geqslant 2$, Alon and R\"{o}dl~\cite{2005AR} proved that $R_{k}(C_{3};K_{m})=\Theta(m^{k+1}\text{poly}\log{m})$ and $R_{k}(C_{2\ell+1};K_m)=\Omega(m^{1+\frac{k}{2\ell-1}}/(\log m)^{k+\frac{2k}{2\ell-1}})$. For more literature on the Ramsey problems, we refer the readers to the dynamic survey~\cite{1994Survey} and the references therein. Inspired by the random homomorphisms method in the recent work~\cite{2021Yuval} on multicolor diagonal Ramsey number, we prove the lower bound for $R_k(H;K_m)$ when $H$ is $C_5$ or $C_7$, respectively. The appropriate base graphs come from the random block construction in~\cite{2019Mubayi}.

\begin{theorem}\label{thm:main}
For fixed positive integer $k$ and $m\rightarrow\infty$, we have
  \begin{equation*}
      R_{k}(C_5;K_m)=\Omega\big((\frac{m}{\log{m}})^{\frac{3k}{8}+1}\big),
  \end{equation*} and
  \begin{equation*}
      R_{k}(C_7;K_m)=\Omega\big((\frac{m}{\log{m}})^{\frac{2k}{9}+1}\big).
  \end{equation*}

\end{theorem}
 These slightly improve the previously known lower bounds $R_{k}(C_5;K_m)=\Omega(m^{\frac{k}{3}+1}\text{poly}\log{m})$ and $ R_{k}(C_7;K_m)=\Omega(m^{\frac{k}{5}+1}\text{poly}\log{m})$ respectively.

\section{Proof of Theorem~\ref{thm:main}}

Our first step is to construct a base graph which is $C_{5}$-free and has no large independent set. The construction is based on a recent work of Mubayi and Verstra\"{e}te~\cite{2019Mubayi}. By the existence of generalized hexagons of order $(q,q^{3})$ (see~\cite{2016IG, 2000JG, 1998book}), there exists a graph $G$ with bipartition $G=A\cup B$ of sizes $|A|=(q+1)(q^{8}+q^{4}+1)$ and $|B|=(q^{3}+1)(q^{8}+q^{4}+1)$ such that every vertex in $A$ has degree $q^{3}+1$ and every vertex in $B$ has degree $q+1$. Moreover, the girth of $G$ is at least $12$, that is, there is no cycle of length less than $12$ in $G$. For a vertex $a\in A,$ we can randomly partition the neighbors $N_{G}(a)$ into two parts $S_{a}$ and $T_{a}$ uniformly, and then we form a complete bipartite graph between $S_{a}$ and $T_{a}$. Furthermore, we repeat these operations for all vertices in $A$ and finally obtain a new graph $F$ on vertex set $V(F)=B.$ Note that there is no multiple edge in the new graph $F$ since $G$ is $C_{4}$-free. Moreover, $F$ is $C_{5}$-free, otherwise by definition of graph $F$, there will be a copy of $C_{10}$ in $G$, a contraction. We also claim that there is no independent set of size $t=(1+o(1))q^{8}$ with positive probability. To see this, let $I$ be a subset of $F$ with $|I|=t$. Set $t_{a}=|I\cap N_{G}(a)|$ for $a\in A,$ it is easy to see the probability that $I\cap N_{G}(a)$ is an independent set in $F$ is $2^{1-t_{a}}$. Since the partitions are independently at random for different $a\in A$, the probability that $I$ forms an independent set in $F$ is $\prod\limits_{a\in A}2^{1-t_{a}}$ and the expected number of independent sets of size $t$ in $F$ is $\binom{|F|}{t}\prod\limits_{a\in A}2^{1-t_{a}}<1.$

Now we fix a graph $F$ which is $C_{5}$-free and has no independent set of size $t=(1+o(1))q^{8}.$ We say $F'$ is an $r$-blowup of the graph $F$ means the graph is obtained by replacing each vertex of $F$ by an independent set of size $r$ and each edge of $F$ by a complete bipartite graph. Set $r=c_{k}q^{3(k-1)}$ and $m=d_{k}q^{8}\log{q}$ for some suitable constants $c_{k},d_{k}>0$. We claim that the number of independent sets of size $m$ in $F'$ is at most $\binom{|F|}{t}\binom{tr}{m}$. To see this, since there is no independent set of size $t$ in $F$, there are at most $\binom{|F|}{t}$ blocks which contain the independent sets in $F'$. Moreover, each vertex of the independent set in $F'$ of size $m$ is one of the $tr$ vertices in these blocks since each block contains exactly $r$ vertices.

We then show that there exists a sequence of $k+1$ graphs $F_{1},F_{2},\ldots,F_{k},F_{k+1}$ on vertex set $V(F')$, such that $F_{i}$ is isomorphic to $F'$ with $1\leqslant i\leqslant k$ and $F_{k+1}$ is the graph whose edges are all pairs of vertices of $V(F')$ that do not lie in any $F_{i}$ with $1\leqslant i\leqslant k$. Moreover, $F_{k+1}$ is $K_{m}$-free. To show this, we pick $k$ random copies $F_{1},F_{2},\ldots, F_{k}$ of $F'$ on $V(F')$ independently, then the graph $F_{k+1}$ is determined by the above rules. For a fixed subset $M$ in $V(F')$ with $|M|=m$, the probability that $M$ forms an independent set in any $F_{i}$ is $\frac{\binom{|F|}{t}\binom{tr}{m}}{\binom{|F'|}{m}}$, thus the probability that $M$ forms a clique in $F_{k+1}$ is at most $\big(\frac{\binom{|F|}{t}\binom{tr}{m}}{\binom{|F'|}{m}}\big)^{k}$, and the expected number of $K_{m}$ in $F_{k+1}$ is at most $\binom{|F'|}{m}\big(\frac{\binom{|F|}{t}\binom{tr}{m}}{\binom{|F'|}{m}}\big)^{k}<1$ by choices of the parameters $r$ and $m$. Hence, the required graph $F_{k+1}$ exists, with positive probability.

Finally, we fix a required sequence of graphs $F_{1},F_{2},\ldots,F_{k+1}$ as above. Let $N:=r|F|$ and color each edge of $K_{N}$ by minimum index $i\in\{1,2,\ldots,k\}$ such that the edge belongs to $F_{i}$. Moreover, if no such $i$ exists, we color the edge by color $k+1$. Through the above analysis, we can see there is no $C_{5}$ in color $i$ with $1\leqslant i\leqslant k$ and no $K_{m}$ in color $k+1$. This gives $R_{k}(C_5;K_m)=\Omega\big((\frac{m}{\log{m}})^{\frac{3k}{8}+1}\big)$.

We can also obtain the lower bound for $R_{k}(C_7;K_m)$ through the similar analysis. The required base graphs of girth at least $16$ can be constructed from Ree-Tits octagons (see~\cite{2016IG, 2000JG, 1998book}), we omit the details.

\section{Remarks}
In~\cite{2021Yuval}, Wigderson pointed out that the random induced subgraphs and random blowups are closely related and are both part of a more general framework of random homomorphisms. We do agree that the framework of random homomorphisms is great, and in some sense we can regard the random induced subgraphs and the random blowups as the same thing. However, in this paper, we prefer to use the random blowups method, since our $C_{5}$-free base graph from the random block construction has too many independent sets of small order. In diagonal Ramsey case~\cite{2021conlon, 2021Yuval}, the base graph is constructed algebraically, which provides a good estimation of the number of small independent sets. The random blowups method has certain limitations, for example, it cannot work when the graph $H$ is bipartite, unless one chooses other ways to define the ``blowups''.

It seems to be very hard to obtain the tight asymptotic bound for $R_{k}(C_{2\ell+1};K_{m})$, even for the simplest case $R_{1}(C_{5};K_{m})$. The best known upper bound $R_{1}(C_{5};K_{m})=O(\frac{m^{3/2}}{\sqrt{\log{m}}})$ was proven in~\cite{2000DM}. For multicolor case, we conjecture that $R_{k}(C_{5};K_{m})=O(m^{\frac{k}{2}+1}\text{poly}\log{m})$ holds.

\section*{Acknowledgements}
The authors thank Xin Wei for carefully reading an early draft of this paper.

\bibliographystyle{abbrv}
\bibliography{Multicolor_Ramsey_cycle}
\end{document}